\documentclass[11pt,twoside]{article}
\usepackage{latexsym,amssymb,amsmath}
\usepackage{amsfonts}
\RequirePackage{graphics,epsfig,psfrag}
\def\abs#1{\left \vert #1 \right \vert}
\def\hat#1{\widehat{#1}}
\def\RR{{\bf R}} 
\def\ZZ{{\bf Z}} 
\def\QQ{{\bf Q}} 
\def\SS{{\bf S}} 
\def\TT{{\rm T}}
\def\GG{{\mathbf{G}}} 
\def\Mod#1{\,(\hbox{\rm mod}\,#1)}
\def\cn{\hbox{\rm cn}\,}
\def\H{{\rm H}}
\def\id{{\mathbf{Id}}}

\def\phi{\varphi}
\def\eps{\varepsilon}

\def\cC{{\cal C}}

\def\cT{{\cal T}}
\def\pn{\medskip\par\noindent}
\def\bi{\vspace{-2pt}\begin{itemize}\itemsep -2pt plus 1pt minus 1pt}
\def\ei{\end{itemize}\vspace{-4pt}}
\def\bn{\vspace{-2pt}\begin{enumerate}\itemsep -2pt plus 1pt minus 1pt}
\def\en{\end{enumerate}\vspace{-4pt}}
\newcommand{\Pf}{{\em Proof}. }

\newcommand{\EPf}{\hbox{}\hfill$\Box$\vspace{.5cm}}
\def\[#1\]{\begin{eqnarray}#1\end{eqnarray}}
\def\$#1\${\begin{eqnarray*}#1\end{eqnarray*}}
\def\p{\vspace{-5.17ex} \hspace{6.5ex}}

\def\sign#1{{\rm sign}\,\bigl( #1 \bigr)}

\def\abs#1{\left \vert #1 \right \vert}

\def\frac#1#2{{\textstyle{{#1} \overwithdelims.. {#2}}}}
\def\Frac#1#2{{\displaystyle{{#1} \overwithdelims.. {#2}}}}

\catcode`\@=11
\def\system#1{\left\{\null\,\vcenter{\openup\jot\m@th
\ialign{
\strut\hfil$\displaystyle{##}$&
$\,\displaystyle{{}##}\,$\hfil&&
\strut\hfil$\,\displaystyle{##}$&
$\,\displaystyle{{}##}\,$\hfill
\hfil\crcr#1\crcr}}\right.}

\def\cmatrix#1{\left [
\null\,\vcenter{
\ialign{
\hfil${##}\ $\hfil &
\hfil$\ {##}\ $\hfil&&
\hfil$\ {##}\ $\hfil&
\hfil$\ {##}$\hfil
\crcr#1\crcr}}\right ]}

\def\@opargbegintheorem#1#2#3{\par\addvspace{6pt plus3pt minus2pt}%
    \def\@tempa{#3}%
    \noindent{\bf #1 #2 \ifx\@tempa\empty\unskip\else\unskip\ (#3).\fi\hskip.5em}\csname#1font\endcsname\ignorespaces
\ignorespaces}

\def\@endtheorem{\par\addvspace{6pt plus3pt minus2pt}}
\def\@begintheorem#1#2#3{\par\addvspace{8pt plus3pt minus2pt}%
              \noindent{\csname#1headfont\endcsname#1\ \ignorespaces#3 #2.}%
              \csname#1font\endcsname\hskip6pt\ignorespaces}
\def\@endtheorem{\par\addvspace{8pt plus3pt minus2pt}\@endparenv}

\catcode`\@=12

\newtheorem{theorem}{Theorem}[section]
\newtheorem{thm*}{Theorem}
\newtheorem{corollary}[theorem]{Corollary}
\newtheorem{algorithm}[theorem]{Algorithm}
\newtheorem{lemma}[theorem]{Lemma}
\newtheorem{proposition}[theorem]{Proposition}
\newtheorem{definition}[theorem]{Definition}
\newtheorem{remark}[theorem]{Remark}
\newtheorem{example}[theorem]{Example}
\newtheorem{examples}[theorem]{Examples}

\newtheorem{conjecture}[theorem]{Conjecture}
\date{\today}
\begin{document}
\pagestyle{myheadings}
\markboth{P. -V. Koseleff, D. Pecker}{{\em Chebyshev diagrams for two-bridge  knots}}
\title{Chebyshev diagrams for two-bridge  knots}
\author{P. -V. Koseleff, D. Pecker\medskip\\
Universit{\'e} Pierre et Marie Curie\\
4, place Jussieu, F-75252 Paris Cedex 05 \\
e-mail: {\tt\{koseleff,pecker\}@math.jussieu.fr}}
\vspace{-1cm}
\maketitle
\begin{abstract}
We show that every two-bridge knot $K$ of crossing number $N$
admits a polynomial parametrization $x=T_3(t), \, y = T_b(t), z =
C(t)$ where $T_k(t)$ are the Chebyshev polynomials and $b+
\deg C = 3N$.
If $C (t)= T_c(t)$ is a Chebyshev polynomial, we call such a knot a harmonic
knot. We give the classification of  harmonic knots for $a \le 3.$
Most results are derived from
continued fractions and their matrix representations.
\pn {\bf keywords:} {Polynomial curves, Chebyshev polynomials, Chebyshev
curves, rational knots, continued fractions}
\end{abstract}
\begin{center}
\parbox{12cm}{\small
\tableofcontents
}
\end{center}
\vspace{1cm}
\section{Introduction}
We study the polynomial parametrization of knots,
viewed as non singular space curves.
Vassiliev proved
that any  knot  can be represented by a polynomial
embedding $\RR \to \RR ^3 \subset {\bf S}_3 $ (\cite{Va}).
Shastri (\cite{Sh}) gave another proof of this theorem,
he also found explicit  parametrizations of the trefoil
and of the figure-eight knot.

We shall study polynomial embeddings of the  form
$x=T_a(t), \, y= T_b(t), \, z=C(t)$ where $a$ and $b$ are coprime
integers and $T_n(t)$ are  the classical Chebyshev polynomials defined by
$T_n(\cos t) = \cos nt$.
The projection of such a curve on the $xy$-plane is
the Chebyshev  curve  $\cC(a,b): T_b(x)=T_a(y)$ which has
exactly $\frac 12 (a-1)(b-1)$ crossing points (\cite{Fi,P1}).
We will say that such a knot has the Chebyshev diagram $\cC(a,b)$.
\pn
We observed in \cite{KP1} that  the trefoil can be parametrized by
Chebyshev polynomials: $x=T_3(t);\, y=T_4(t);\,z= T_5(t)$.
This led us to study Chebyshev knots in \cite{KP3}. We obtained the following result:
\begin{theorem}[\cite{KP3}]
Any knot is a Chebyshev knot, that is, is isotopic to a knot given
by a one-to-one parametrization
$$ \cC(a,b,c,\phi): \ x=T_a(t); \  y=T_b(t) ; \    z=T_c(t + \phi) $$
where $t \in \RR$, $a$ and $b$ are coprime integers, $c$
is an integer and $\phi$ is a real constant.
\end{theorem}
Our proof uses theorems on braids by Hoste, Zirbel and Lamm (\cite{HZ,La2}),
and a density argument.
In a joint work with F. Rouillier (\cite{KPR}), we developed an effective method
to enumerate all the knots $\cC(a,b,c,\phi), \phi \in \RR$ where $a=3$ or $a=4$, $a$ and
$b$ coprime.
\pn
Chebyshev knots are polynomial analogues  of
Lissajous knots that admit a  parametrization of the form
$$
x=\cos (at ); \   y=\cos (bt + \phi) ; \  z=\cos (ct + \psi)
$$
where $ 0 \le t \le 2 \pi $ and where $ a, b, c$ are pairwise
coprime integers. These knots, introduced  in \cite{BHJS},
have been studied by  V. F. R. Jones, J. Przytycki,
C. Lamm, J. Hoste and L. Zirbel. Most known properties of Lissajous
 knots are deduced from their symmetries
(see \cite{BDHZ,Cr,HZ,JP,La1}).
\pn
\begin{definition}
When $a,b,c$ are coprime then $\cC(a,b,c,0)$ is
denoted by $\H(a,b,c)$ and is called a harmonic knot.
\end{definition}
The symmetries of harmonic knots, obvious from the parity of Chebyshev polynomials,
are different from those of Lissajous.
For example, the figure-eight knot which is amphicheiral but not a Lissajous knot, is
the harmonic knot $ \H(3,5,7).$
\pn
We proved in \cite{KP3} that the harmonic knot $\H(a,b, ab-a-b)$ is alternating, and deduced
that there  are infinitely many amphicheiral harmonic knots and
infinitely many strongly invertible harmonic knots.
We also proved in \cite{KP3} that the torus knot $\TT(2, 2n+1)$ is the harmonic knot
$ \H(3,3n+2,3n+1)$.
\pn
In this article, we give the classification of the harmonic knots $\H(a,b,c)$ for $a \leq 3.$
We also give explicit polynomial parametrizations of all rational knots.
The diagrams of our knots are Chebyshev curves of minimal degrees with a small number
of crossing points. The degrees of the height polynomials are small.
\pn
In section {\bf \ref{cf}} we recall the Conway notation for rational knots,
and the computation of their Schubert fractions with continued fractions.
We observe that, when $a=3$, Chebyshev diagrams correspond
to continued fractions of the form $[\pm 1, \ldots, \pm 1]$.
\pn
The study of these particular continued fraction expansions will be the main tool
of this article.
\pn
{\bf Theorem \ref{th1}.}\\
{\em Every rational number $r$ has a unique continued fraction expansion
$r = [e_1, e_2, \ldots, e_n ]$, $e_i= \pm 1$,
where there are no two consecutive sign changes in the sequence $(e_1,\ldots , e_n ).$ }
\pn
We provide a formula (Proposition \ref{bireg}) for the crossing number of the corresponding knots.
Then we study the matrix interpretation of these continued fraction expansions.
\pn
In section {\bf \ref{diagrams}} we show how to find explicit
minimal Chebyshev diagrams $\cC(3,b)$ for all rational knots:
\pn
{\bf Theorem \ref{ell}.}\\
{\em Let $K$ be a two-bridge knot with crossing number $N$.
There is an algorithm to determine the smallest $b$ such that
$K$ has a Chebyshev diagram $\cC(3,b)$ with $N < b < \Frac 32 N$.}
\pn
As an application, we give optimal Chebyshev diagrams for the torus knots
$\TT(2,N)$, the twist knots $\cT_n$, the generalized stevedore knots and some
others.
\pn
In section {\bf \ref{param}}, we find explicit polynomial parametrizations of all
rational knots. We first use the minimal Chebyshev diagram found in section {\bf \ref{diagrams}},
then we define a height polynomial of small degree.
More precisely, we show:
\pn
{\bf Theorem \ref{gauss3}.}\\
{\em Every rational knot of crossing number $N$ can be parametrized by
$x=T_3(t), y= T_b(t), z= C(t)$ where $b + \deg C = 3N$.
Furthermore,  when the knot
is amphicheiral, $b$ is odd and we can choose $C$ to be an odd polynomial.}
\pn
We give the first polynomial parametrizations of the twist knots and the generalized stevedore knots.
\pn
In section {\bf \ref{harmonic}}, we give a complete classification of
harmonic knots $\H(3,b,c)$ where $b$ and $c$ are relatively prime integers, not divisible
by 3. We obtain this classification by a careful study of the diagrams and
their continued fractions  (Theorem \ref{h3} whose proof is given in section {\bf \ref{proofs}}).
We show that the twist knots and the generalized stevedore knots (e.g. the $6_1$ knot)
are not harmonic knots $\H(3,b,c)$.
\pn
Thanks to the use of continued fractions, we provide effective methods for the construction
of polynomial parametrizations for any rational knot. We conjecture they are of minimal
degrees.
\section{Continued fractions and rational Chebyshev knots}\label{cf}
\begin{figure}[th]
\begin{center}
\begin{tabular}{cccc}
{\scalebox{.6}{\includegraphics{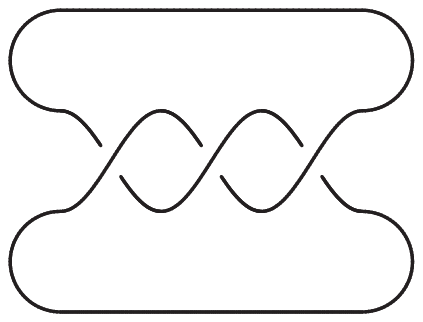}}}&
{\scalebox{.6}{\includegraphics{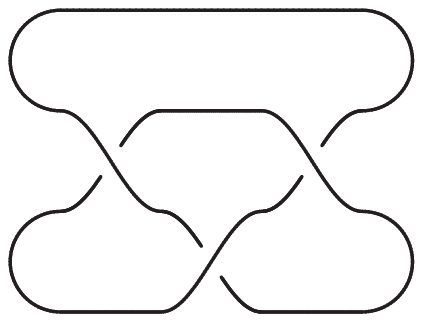}}}&
{\scalebox{.6}{\includegraphics{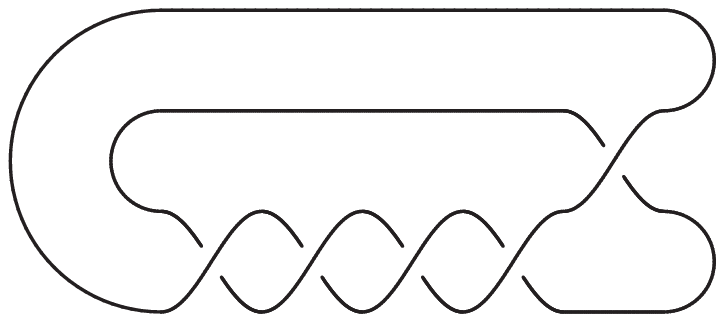}}}&
{\scalebox{.6}{\includegraphics{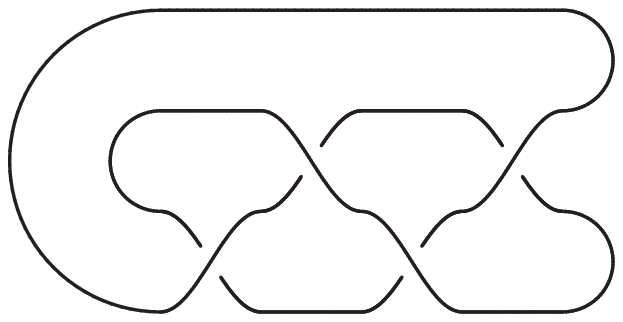}}}\\
$C(3)$&$C(-1,-1,-1)$&$C(4,-1)$&$C(1,1,-1,-1)$
\end{tabular}
\end{center}
\caption{Diagrams of the standard trefoil}\label{trefoils}
\end{figure}

A two-bridge knot (or link) admits a diagram in Conway's normal form.
This form, denoted by
$C(a_1, a_2,\ldots, a_n)$  where $a_i$ are integers, is explained by the following
picture (see \cite{Con}, \cite{Mu} p. 187).
\psfrag{a}{\small $a_1$}\psfrag{b}{\small $a_2$}%
\psfrag{c}{\small $a_{n-1}$}\psfrag{d}{\small $a_{n}$}%
\begin{figure}[th]
\begin{center}
{\scalebox{.8}{\includegraphics{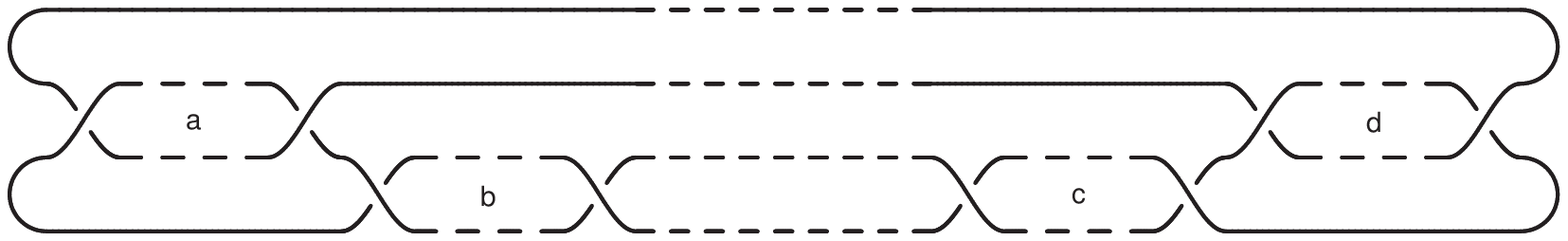}}}\\[30pt]
{\scalebox{.8}{\includegraphics{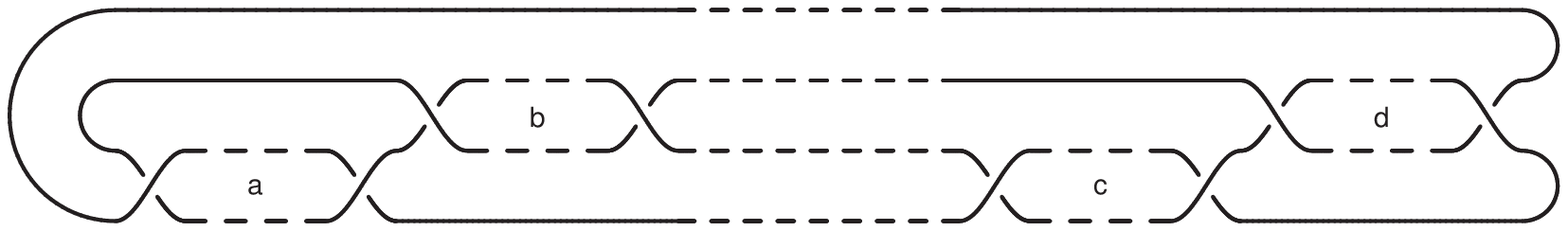}}}
\end{center}
\caption{Conway's normal forms, $n$ odd, $n$ even}
\label{conways3}
\end{figure}
The number of twists is denoted by the integer
$\abs{a_i}$, and the sign of $a_i$ is defined
as follows: if $i$ is odd, then the right twist is positive,
if $i$ is even, then the right twist is negative.
On  Fig. \ref{conways3} the $a_i$ are positive (the $a_1$ first twists are right twists).
\pn
\begin{examples}
The trefoil has the following  Conway's normal forms
$C(3)$, $C(-1,-1,-1)$, $C(4, -1)$ and $C( -1,-1,1,1 ).$
The diagrams in Figure \ref{trefoils} clearly represent the same trefoil.
\end{examples}
\def\p{{\small $+$}}
\def\m{{\small $-$}}
\pn
The two-bridge links are classified by their Schubert fractions
$$
 \Frac {\alpha}{\beta} =
a_1 + \Frac{1} {a_2 + \Frac {1} {a_3 + \Frac{1} {\cdots +\Frac 1{a_n}}}}=
[ a_1, \ldots , a_n] , \quad \alpha >0.
$$
We shall denote  by $S \bigl( \Frac {\alpha}{\beta} \bigr)$  a two-bridge link with
Schubert fraction $ \Frac {\alpha}{\beta} .$
The two-bridge  links
$ S (\Frac {\alpha} {\beta} )$ and $ S( \Frac {\alpha ' }{\beta '} )$ are equivalent
if and only if $ \alpha = \alpha' $ and $ \beta' \equiv \beta ^{\pm 1} ( {\rm mod}  \  \alpha).$
The integer $ \alpha$ is odd for a knot, and even for a two-component link.
If $K= S (\Frac {\alpha}{\beta} ),$ its mirror image is
$ \overline{K}= S ( \Frac {\alpha}{- \beta} ).$

We shall study knots with a Chebyshev diagram $\cC (3,b) : \  x= T_3(t), y= T_b(t).$
It is remarkable that such a diagram is already in Conway normal form: the crossing points belong
to the 2 horizontal lines $y=\pm \frac12 \sqrt3$. 
Consequently, the Schubert fraction of such a knot is given by a
continued fraction of the form
$ [ \pm 1, \pm 1, \ldots ,\pm 1 ] .$
For example the only  diagrams of Figure \ref{trefoils}
which may be Chebyshev are the second and the last
(in fact they are Chebyshev).
\pn
Figure \ref{t7} shows  a typical example of  a knot with a Chebyshev diagram.
\def\p{{$\scriptstyle{+}$}}%
\def\m{{\small $-$}}%
\begin{figure}[th]
\begin{center}
\psfrag{n1}{\m}\psfrag{n2}{\m}\psfrag{n3}{\m}%
\psfrag{n4}{\p}\psfrag{n5}{\p}\psfrag{n6}{\p}%
\psfrag{n7}{\m}\psfrag{n8}{\m}\psfrag{n9}{\m}%
{\scalebox{1}{\includegraphics{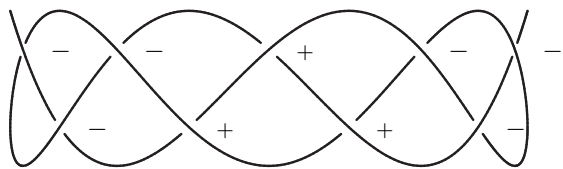}}}
\end{center}
\caption{A Chebyshev diagram of  the torus knot $\TT(2,7)$}
\label{t7}
\end{figure}
\pn
This knot is defined by $ x= T_3(t), \, y= T_{10}(t), \, z = -T_{11}(t) .$
Its $xy$-projection is in the Conway normal form
$ C(-1,-1,-1 , 1,1,1, -1,-1,-1).$ Its Schubert fraction is then $\Frac{7}{-6}$ and this knot is
the torus knot $\TT(2, 7)=S(\Frac{7}{-6}) = S(7)$.
\pn
Let $\alpha,\beta$ be relatively prime integers.
Then $ \Frac {\alpha }{\beta}$ admits the continued fraction expansion
$\Frac{\alpha}{\beta} = [q_1, q_2, \ldots, q_n ] $
if and only if there exist integers  $r_i \in \ZZ$ such that
$$
\left \{
\begin{array}{rcl}
\alpha &=& q_1 \beta +r_2,\\
\beta &=& q_2 r_2 + r_3,\\
&\vdots&\\
r_{n-2} &=& q_{n-1} r_{n-1} + r_n,\\
r_{n-1} &=& q_n r_n.
\end{array}
\right .
$$
The integers $q_i \in \ZZ$ are called the quotients of the continued fraction.
Euclidean algorithms provide various continued
fraction expansions which are useful to the study of two-bridge knots
(see \cite{BZ,St,Cr}).
\pn
If $\alpha>\beta>0$, there is a unique continued expansion of
$\Frac\alpha\beta = [q_1, \ldots, q_n], \, q_i >0$ up to $q_n = [q_{n}-1,1]$.
\begin{definition}
Let $r >1$ be a rational number, and $r=[q_1, \ldots , q_n ] $ be its classical
continued fraction expansion (with $q_i >0$).
The {\em crossing number} of $r$ is defined by
$\cn(r)= q_1 + \cdots + q_n. $
\end{definition}
\begin{remark}\label{cns}
One can prove that $\cn(\Frac\alpha\beta)=\cn(\Frac\alpha{\alpha-\beta})$ when $\Frac\alpha\beta>1$.
If $K= S\bigl ( \Frac\alpha\beta  \bigr)$, $\Frac\alpha\beta>1$ is a rational knot, it
is known that $\cn(\Frac\alpha\beta)$ is the crossing number of $K$.
It means that it is  the minimum number of crossing points for all diagrams of $K$ (\cite{Mu}).
\end{remark}
We shall be interested in algorithms
where the sequence of remainders is not necessarily decreasing anymore (the $q_i$
are not necessarily positive).
\begin{definition}
A continued fraction $ [ a_1,a_2 , \ldots, a_n]$ is 1-regular if it
has the following
properties:
$$
a_i \neq 0, \,  a_{n-1} a_n >0 ,   \hbox{ and }
a_i a_{i+1} <0 \Rightarrow  a_{i+1} a_{i+2} >0, \ i = 1, \ldots, n-2.$$
\end{definition}
\begin{proposition}\label{bireg}
Let $\Frac \alpha\beta = [a_1, \ldots, a_n]$
be a 1-regular continued fraction with $a_1,a_2>0$. Then $\Frac \alpha\beta>1$ and we have
\[
\cn(\Frac \alpha\beta) = \sum_{k=1}^n \abs{a_i} - \sharp \{i,
a_ia_{i+1}<0\}.
\label{biregf}
\]
\end{proposition}
\Pf We prove
this result by induction on the number of sign changes
$k = \sharp \{i, a_ia_{i+1}<0\}$.
If $k$ is 0, it is the definition. If
$k>0$ let us consider the first change of sign. The 1-regular continued fraction is
$[x,a,b,-c,-d,-y]$ where $a,b,c,d$ are positive integers, $x$
is a sequence (possibly empty) of positive integers and $y$ is a 1-regular sequence
of integers.
We have $ [x,a,b,-c,-d,-y] = [x,a,b-1,1,c-1,d,y]$ (Lagrange identity, \cite{Cr})
\bi
\item Suppose $(b-1)(c-1)>0$, then
the sum of absolute values has decreased by 1 and the number of changes of sign
has also decreased by 1.
\item Suppose $b=1, c \not = 1$ (resp. $c=1, b\not = 1$). Then
$[x,a,b,-c,-d,-y] = [x,a,0,1,c-1,d,y] =[x,a+1,c-1,d,y]$. (resp. $[x,a-1,c+1,d,y]$).
The sum of absolute values has decreased by 1 and the number of changes of sign
has also decreased by 1.
\item Suppose $b=c=1$. Then
$[x,a,b,-c,-d,-y] = [x,a,0,1,0,d,y] =[x,a+d+1,y]$.
The sum of absolute values has decreased by 1 and the number of changes of sign
has also decreased by 1.
\ei
We therefore deduce that $\Frac\alpha\beta = [x,r]$ where $r>1$ and $x$ is a sequence
sequence (possibly empty) of positive integers. We get $\Frac\alpha\beta>1$.
\EPf
\pn
Note that Formula (\ref{biregf}) still holds
when $a_1, \ldots ,a_n$ are non zero even integers and the sequence is not
necessarily 1-regular (see \cite{St}).
\pn
We shall now use the basic (subtractive) Euclidean algorithm to get 1-regular continued fractions
of the form $ [ \pm 1, \pm 1 , \ldots , \pm 1 ] .$
\section{Continued fractions ${\mathbf{[\pm 1,\pm 1,\ldots,\pm 1]}}$}\label{pm1}
We will consider the following M\"{o}bius transformations:
\[
P : x \mapsto [1,x] = 1+\Frac 1x, \,
M : x \mapsto [1,-1,-x] = \Frac{1}{1+x}. \label{PM}
\]
Let $E$ be the set of positive real numbers.
We have  $P(E) = ]1, \infty [$ and $M(E)= ]0, 1 [.$
$P(E)$ and $M(E)$ are disjoint subsets of $E.$
\begin{theorem}\label{th1}
Let $ \Frac \alpha\beta$ be a rational number.
There is a unique 1-regular continued fraction such that
$$\Frac \alpha\beta = [e_1, e_2, \ldots,  e_n], \,  e_i= \pm 1.$$
Furthermore, $\Frac \alpha\beta>1$ if and only if $e_1=e_2=1$.
\end{theorem}
\Pf
Let us prove the existence by induction on the height
$h(\Frac\alpha\beta)=\alpha+\beta$.
\bi
\item
If $h=2$ then  $\Frac \alpha \beta = 1 = [1]$ and the  result is true.
\item
If $ \alpha > \beta,$ we have
$\Frac \alpha \beta = P(\Frac \beta {\alpha - \beta}) =
[1 ,\Frac \beta {\alpha - \beta }].
$
Since $h( \Frac \beta {\alpha - \beta } ) < h ( \Frac \alpha \beta )$,
we get our 1-regular continued fraction for $r$ by induction.
\item
If $ \beta > \alpha$  we have
$ \Frac \alpha \beta = M( \Frac  {\beta - \alpha } \alpha ) =
[1, -1, -\Frac {\beta - \alpha} \alpha  ].
$
And we also get a 1-regular continued fraction for $r$.
\ei
This completes
the construction of our continued fraction expansion.
\pn
On the other hand, let $r$ be defined by the 1-regular continued fraction
$ r = [ 1, r_2 , \ldots, r_n]$, $r_i = \pm 1$,  $n \ge 2 $.
Let us prove, by induction on the length $n$ of the continued fraction, that $r>0$
and that $r >1$ if and only if $r_2 = 1$.
\bi
\item If $r_2=1$ we have $r=P([1,r_3,\ldots ,r_n] ),$ and by induction $r\in P(E)$ and
then $r>1$.
\item If $r_2=-1,$ we have $r_3=-1$ and $r= M([1, -r_4,\ldots, -r_n])$.
By induction, $r\in M(E) = ]0,1[$.
\ei
The uniqueness is now easy to prove.
Let $ r = [1, r_2, \ldots , r_n] = [1, r'_2, \ldots , r'_{n'} ].$
\bi
\item If $r>1$ then $r_2=r'_2=1$ and
$[1,1, r_3, \ldots , r_n ] =[1, 1, r'_3, \ldots r'_{n'} ]$.
Consequently,
$[1, r_3, \ldots , r_n] = [1 , r'_3, \ldots, r'_{n'} ],$ and by induction
$r_i=r'_i$ for all $i.$
\item If $r < 1,$ then $r_2=r_3=r'_2=r'_3=-1$ and
$[1,-1,-1, r_4, \ldots , r_n] =[1,-1,-1, r'_4 ,\ldots, r'_{n'} ].$
Then,
$ [ 1, -r_4, \ldots, -r_n] = [ 1, -r'_4, \ldots, -r'_{n'} ]$
and by induction
$ r_i=r'_i$ for all $i.$
\EPf
\ei
\begin{remark}
Since we have $[1,-1,1,x] = -x$, we see that 1-regularity is required
for the uniqueness of our expansion. This is exactly analogous to the classical case.
\end{remark}
\begin{definition}
Let $\Frac\alpha\beta>0$ be the 1-regular continued fraction $[e_1, \ldots, e_n]$, $e_i=\pm 1$.
We will denote its length $n$ by $\ell(\Frac\alpha\beta)$. Note that
$\ell(-\Frac\alpha\beta) = \ell(\Frac\alpha\beta)$.
\end{definition}
\begin{examples}\label{9297}
Using our algorithm we obtain
$$
\begin{array}{rcl}
\Frac 97&=& [1,\Frac 72] = [1,1,\Frac 25] =  [1,1,1,-1, -\Frac 32] =
[1,1,1,-1, -1, -\Frac 21] = [1,1,1,-1, -1, -1, -1],\\[10pt]
\Frac 92&=& [1,\Frac 27] = [1,1,-1,-\Frac 52] = [1,1,-1,-1,-\Frac 23] =
[1,1,-1,-1,-1,1,\Frac 12] \\
&=& [1,1,-1,-1,-1,1,1,-1,-1] = [4,2].
\end{array}
$$
\begin{figure}[th]
\begin{center}
\begin{tabular}{ccc}
{\scalebox{.6}{\includegraphics{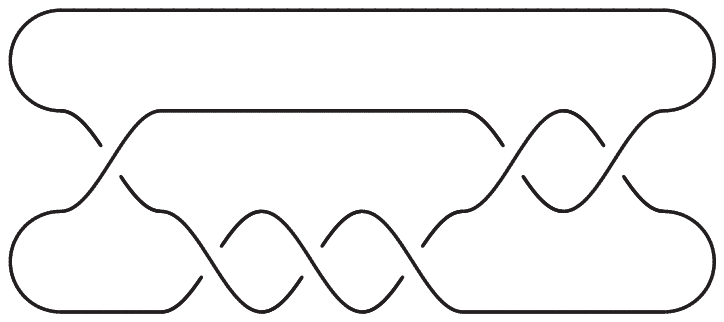}}}&
{\scalebox{.6}{\includegraphics{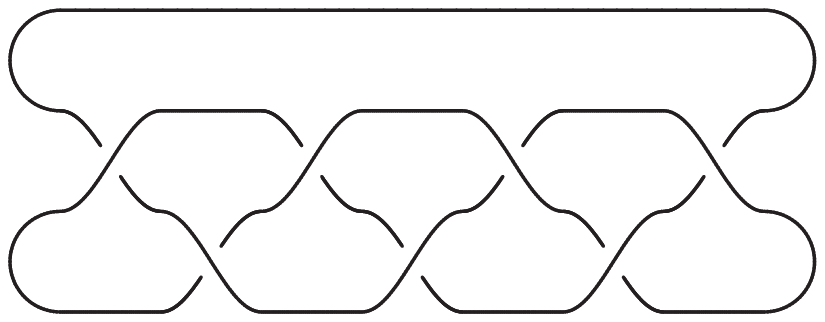}}}&
{\scalebox{.6}{\includegraphics{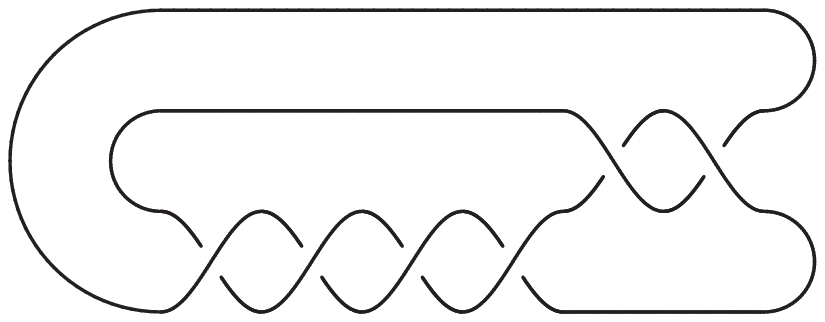}}}\\
$C(1,3,2)$&$C(1,1,1,-1,-1,-1,-1)$&$C(4,2)$
\end{tabular}
\end{center}
\caption{Diagrams of the knot $6_1 = S(\Frac 97)$ and its mirror image $S(\Frac 92)$ }
\label{k61}
\end{figure}
\end{examples}
We will rather use the notation
$$
\Frac 97 = P^2 M P^3 (\infty), \
\Frac 92 = P MP M^2 P (\infty).
$$
We get $\ell (\Frac 97) = 7$, $\ell (\Frac 92) = 9$.
The crossing numbers of these fractions are
$\cn(\Frac 97) = \cn([1,3,2]) = 6 = 7-1$ and $\cn(\Frac 92) = \cn([4,2]) = 6 = 9-3$.
If the fractions $\Frac 97$ and $\Frac 92$ have the same crossing number,
it is because the knot $S(\Frac 97)$ is the mirror image of $S(\Frac 92)$.
In order to get a full description of two-bridge knots we shall need a more detailed
study of the homographies $P$ and $M$.
\pn
\begin{proposition}
The multiplicative monoid $\GG=\langle P,M \rangle$
is free.
The mapping $g: G \mapsto G(\infty)$ is a bijection from $\GG \cdot P$ to $\QQ_{>0}$ and
$g(P\cdot \GG \cdot P) = \QQ_{>1}$, the set of rational  numbers greater than $1.$
\end{proposition}
\Pf Suppose  that $PX=MX'$ for some $X,X'$ in $\GG$. Then we
would have $PX(1)=MX'(1) \in P(E)\bigcap M(E) = \emptyset$. Clearly, this
implies that $\GG$ is free. Similarly, from $P(\infty)=1$, we deduce that the
mapping $G \mapsto G \cdot P (\infty)$ is injective. From
Theorem \ref{th1} and $P(\infty)=1$,  we deduce that $g$ is
surjective. \EPf
\begin{remark}\label{rl}
Let $r=G(\infty)=[e_1, \ldots, e_n] , \  e_i=\pm1 , $ be a 1-regular continued fraction.
It is easy to find the unique homography $G \in \GG\cdot P$ such
that $r=G(\infty)$. Consider the sequence $(e_1, \ldots, e_n)$. For any
$i$ such that $e_ie_{i+1}<0$, replace the couple $(e_i,e_{i+1})$ by
$M,$ and then replace each  remaining $e_i$ by $P$.
\pn
Let $G= P^{p_1} M^{m_1} \cdots M^{m_k} P^{p_{k+1}}$.
Let $p= p_1 + \cdots + p_{k+1}$ be the degree
of $G$ in $P$ and $m = m_1 + \cdots + m_k$ its degree in $M$. Then we have
$n = \ell(r) = p + 2m$ and $\cn(r) = p+m$.
\end{remark}
We shall consider matrix notations for many proofs. We will consider
$$\cmatrix{\alpha\cr\beta} = P^{p_1} M^{m_1} \cdots M^{m_k} P^{p_{k+1}}
\cmatrix{1\cr 0}, \quad
P = \cmatrix{1&1\cr1&0},\quad M =\cmatrix{0&1\cr1&1}.
$$
\begin{definition}\\
We define on $\GG$ the anti-homomorphism $G\mapsto \overline G$ by $\overline M=M$, $\overline P=P$. \\
We define on $\GG$ the homomorphism $G \mapsto \hat G$ by $\hat M = P, \, \hat P = M$.
\end{definition}
\begin{proposition}\label{gg}
Let $\alpha>\beta>0$ and consider $\Frac\alpha\beta=PGP(\infty)$ and $N=\cn(\Frac\alpha\beta)$.
Let $\beta'$ be such that $0<\beta'<\alpha$ and $\beta\beta'\equiv (-1)^{N-1}\Mod\alpha$.
Then we have
$$
\Frac\beta\alpha = M\hat G P(\infty),\quad
\Frac\alpha{\alpha-\beta} = P \hat G P(\infty), \quad
\Frac\alpha{\beta'} = P \overline G P(\infty).$$
We also have
$$
\ell (\Frac {\beta}{\alpha}) + \ell( \Frac {\alpha}{\beta }) = 3 N -1, \quad
\ell ( \Frac {\alpha}{\alpha - \beta } ) + \ell ( \Frac {\alpha}{\beta})
= 3 N -2,
\quad
\ell ( \Frac {\alpha}{\beta'} ) = \ell( \Frac {\alpha}{\beta }).
$$
\end{proposition}
\Pf
We use matrix notations for this proof.
Let us consider
$PGP=\cmatrix{\alpha&\beta'\cr\beta&\alpha'}= P^{p_1} M^{m_1} \cdots M^{m_k} P^{p_{k+1}}$.
From $\det P = \det M  = -1$, we obtain
$\alpha\alpha' - \beta\beta' = (-1)^N$.
Let $A =\cmatrix{a&c\cr b&d}$ be a matrix such that
$0\leq c \leq a$ and $0\leq d \leq b$.
From $PA = \cmatrix{a+c&b+d\cr a&c}$ and
$MA = \cmatrix{b&d\cr a+b&c+d}$ we deduce that
$PGP$ satisfies $0<\alpha'<\beta$ and $0<\beta'<\alpha$.
We therefore conclude that,
$\beta'$ is the integer defined by
$0<\beta'<\alpha$, $\beta\beta'\equiv (-1)^{N-1} \Mod\alpha$.
By transposition we deduce that
$$\cmatrix{\alpha&\beta\cr\beta'&\alpha'} =
P^{p_{k+1}} M^{m_{k}} \cdots M^{m_1} P^{p_1} = P \overline G P,$$
which implies $\Frac{\alpha}{\beta'} = P \overline G P(\infty)$.
\pn
Let us introduce $J = \cmatrix{0&1\cr1&0}$. We have $J^2 = \id, M=JPJ$ and $P=JMJ.$
Therefore
$$
\cmatrix{\beta\cr\alpha} = J \cmatrix{\alpha\cr\beta}=
M^{p_1} P^{m_1} \cdots P^{m_k} M^{p_{k+1}-1} J P \cmatrix{1 \cr 0}=
M^{p_1} P^{m_1} \cdots P^{m_k} M^{p_{k+1}-1} P \cmatrix{1 \cr 0},
$$
that is $\Frac\beta\alpha = M\hat G P (\infty) $.
\pn
Then,
$\cmatrix{\alpha \cr \alpha - \beta } = PM^{-1} \cmatrix{ \beta \cr \alpha }
= P \hat G P \cmatrix{1 \cr 0} .$
That is
$ \Frac {\alpha}{ \alpha - \beta} = P \hat G P ( \infty ) .$
\pn
Relations on lengths are derived from the previous relations and remark \ref{rl}.
\EPf
\begin{remark}\label{cns2}
It is straightforward that if $[\eps_1, \ldots, \eps_n] = PGP(\infty)$ then
$P\overline G P(\infty) = \eps_n [\eps_n, \ldots, \eps_1]$. We deduce from
Proposition \ref{bireg} that $\cn(\Frac {\alpha}{\beta})=\cn(\Frac {\alpha}{\beta'})$, as expected.
\end{remark}
\begin{lemma}\label{length}
Let $\Frac{\alpha}{\beta} = [e_1, \ldots, e_n]$ be a 1-regular
continued
fraction ($e_i = \pm 1$). We have
\bi
\item $n\equiv 2 \Mod 3$ if and only if $\alpha$ is even and $\beta$ is odd.
\item $n\equiv 0 \Mod 3$ if and only if $\alpha$ is odd and $\beta$ is even.
\item $n\equiv 1 \Mod 3$ if and only if $\alpha$ and $\beta$ are odd.
\ei
\end{lemma}
\Pf
Let us write  $\cmatrix{\alpha\cr\beta} = P^{p_1} M^{m_1} \cdots M^{m_k} P^{p_{k+1}}
\cmatrix{1\cr 0}$.
Since $M \equiv P^2 \Mod 2$, and $n= p+2m$, we get
$P^{p_1} M^{m_1} \cdots M^{m_k} P^{p_{k+1}} \equiv P^n \Mod 2$.
As $P^3 \equiv \id \Mod 2 $ we obtain
\bi
\item[]
if $n\equiv 2 \Mod 3$,
then $\cmatrix{\alpha\cr \beta} \equiv M \cmatrix{1\cr 0} \equiv \cmatrix{0\cr 1} \Mod 2$,
\item[]
if $n\equiv 1 \Mod 3$, then $\cmatrix{\alpha\cr \beta} \equiv P \cmatrix{1\cr 0} \equiv \cmatrix{1\cr 1} \Mod 2$,
\item[]
if $n\equiv 0 \Mod 3$, then $\cmatrix{\alpha\cr \beta} \equiv \cmatrix{1\cr 0} \Mod 2$.
\EPf
\ei
We deduce the following useful result
\begin{proposition}\label{palin}
Let $G \in \GG $ and  $\Frac\alpha\beta = [e_1, \ldots, e_n] = PGP(\infty).$
Let $K = S( \Frac \alpha\beta )$ and $N= \cn (\Frac \alpha\beta)$.
The following properties are equivalent:
\bn
\item $G$ is palindromic (i.e. $\overline G = G$).
\item the sequence of sign changes in $[e_1, \ldots, e_n]$ is palindromic
(i.e. $e_ie_{i+1} = e_{n-i}e_{n-i+1}$).
\item $\beta^2 \equiv (-1)^{N-1}\Mod \alpha$.
\en
Furthermore we have
\bi
\item
$\beta^2 \equiv -1 \Mod \alpha$ (i.e. $K= \overline{K}$ is amphicheiral)
if and only if $N$ is even and $G = \overline G$. Furthermore, the length
$n=\ell( \frac \alpha \beta) $ is even and the sequence $[e_1, \ldots, e_n]$
is palindromic (i.e. $e_i = e_{n-i+1}$).
\item
$\beta^2\equiv 1 \Mod\alpha$ if and only if $N$ is odd and $G = \overline G$
or $N$ is even and $\hat G = {\overline G}$ (in this case $K$ is a 2-component link).
\ei
\end{proposition}
\Pf
From Remark \ref{cns2}, we deduce that
$G=\overline G$ is palindromic if and and only if the sequence of sign changes
in $[e_1, \ldots, e_n]$ is palindromic.

Let $0<\beta'<\alpha$ such that $\beta'\beta \equiv (-1)^{N-1}\Mod\alpha$.
We have from the previous proposition: $\Frac\alpha{\beta'} = P \overline{G}P ( \infty )$.
We thus deduce that  $G=\overline G$ is equivalent to $\beta = \beta'$, that is
$\beta^2\equiv(-1)^{N-1}\Mod\alpha$.

Suppose now that $\beta^2\equiv 1\Mod\alpha$. If $N$ is even then $\beta'=\alpha-\beta$,
that is $P\overline GP = P\hat G P$ and $\overline G=\hat G$.
We  have $p+2m = m + 2 p-2$ and then $2n =2(p+2m) = 3N-2$. This implies
$n\equiv 2\Mod 3$. By  Lemma \ref{length}, $\alpha$ is even and $K$ is a two-component link.
If $N$ is odd then $\beta'=\beta$ and $G=\overline G$ by the first part of our proof.

Suppose now that $\beta^2\equiv -1 \Mod\alpha$.
If $N$ is odd then $\beta'=\alpha-\beta$ and
by the same argument we should have $n=3N-n-2,$
which would imply that $N$ is even. We deduce
that amphicheiral rational links have even crossing numbers
and from $\beta'=\beta$ we get $G=\overline G$.
The crossing number $N=m+p$ is even and $G$ is palindromic so $m$ and $p$ are both even.
Consequently $n= p+2m$ is even and the number of sign changes is even. We thus have
$e_n=1$ and $(e_n, \ldots, e_1)=(e_1, \ldots, e_n)$, using remark \ref{cns2}.
\EPf
\section{Chebyshev diagrams of rational  knots}\label{diagrams}
\begin{definition} We say that a knot in $\RR ^3  \subset  {\SS}^3$ has
a Chebyshev diagram $\cC(a,b)$, if $a$ and $b$ are coprime and the Chebyshev curve
$$\cC(a,b): x=T_a(t); \  y=T_b(t)$$ is the projection of some knot which is
isotopic to $K$.
\end{definition}
In \cite{KP3} we proved that every knot has a Chebyshev diagram.
\begin{proposition}
Let $K$ be a knot, $br(K)$ its bridge number. Let $m \geq br(K)$ be an integer.
Then $K$ has a projection which is a Chebyshev curve $\cC(a,b): x=T_a(t); \  y=T_b(t)$,
where $a = 2m-1$ and $b \equiv  2 \Mod{2a}$.
\end{proposition}
This result is analogous to a theorem of Lamm for
Lissajous curves (see \cite{La2,BDHZ}).
\pn
In the case of two-bridge knots, we give an easy proof of this result.
Moreover, we give an explicit method to get a minimal Chebyshev diagram $\cC(3,b)$.
\begin{theorem}\label{ell}
Let $K$ be a two-bridge knot with crossing number $N$.
\bn
\item $K$ has a Chebyshev diagram $\cC(3,b)$ with $N < b < \Frac 32 N$.
\item There exists $\Frac\alpha\beta>1$ such that $K=S(\pm\Frac\alpha\beta)$ and
$\ell(\Frac\alpha\beta)<\frac 32 N-1$.
If $\Frac\alpha\beta$ is such a fraction, then
$b=\ell(\Frac\alpha\beta)+1$ is the minimal
integer such that $K$ has a Chebyshev diagram $\cC(3,b)$.
\item If $K$ has the Conway normal forms $C(\eps_1,\ldots,\eps_n)$, $\eps_i=\pm 1$,
and $C(e_1,\ldots,e_n)$, $e_i =\pm 1$ of minimal length,
then we have either
$(e_1, \ldots, e_n)=(\eps_1,\ldots,\eps_n)$ or
$(e_1, \ldots, e_n)=(-1)^{n+N}(\eps_n,\ldots,\eps_1)$.
\en
\end{theorem}
\Pf
Let $K$ be a two-bridge knot. Let $r=\Frac\alpha\beta>1$ such that $K=S(r)$.
We have $\overline K = S(r')$ where $r'=\Frac{\alpha}{\alpha-\beta}$.
From Proposition \ref{gg} and Proposition \ref{bireg}, we have
$\ell(r)+\ell(r')=3N-2$ and therefore $N\leq \min(\ell(r),\ell(r')) < \frac 32 N$.
From Lemma \ref{length}, we have $\ell(r)\not \equiv 2 \Mod 3$ so $\ell(r)\not=\ell(r')$
and $n=\min(\ell(r),\ell(r'))<\frac 32 N -1$.
Let us suppose now that $n=\ell(r)<\ell(r')$ and consider the 1-regular continued fraction expansion
$r=[e_1, \ldots, e_n]$. Then $C(e_1,\ldots,e_n)$ is a Conway normal
form for $K$. This Conway normal form corresponds to a Chebyshev diagram
$\cC(3,n+1): x= T_3(t),\, y=T_{n+1}(t)$ and $b=n+1$. If $n=\ell(r')$ and
$r'=[e'_1, \ldots, e'_n]$, we would have considered the Conway normal
form $C(-e'_1,\ldots,-e'_n)$ for $K=S(-r')$.
\pn
Let us consider $\gamma$ such that
$\beta\gamma \equiv 1 \Mod\alpha$ and $0<\gamma<\alpha$.
Let $\rho = \Frac{\alpha}{\gamma}$ and $\rho' = \Frac{\alpha}{\alpha-\gamma}$.
We have $K = S(\rho)$ and $\overline K = S(\rho')$ and from Proposition \ref{gg}:
$n=\ell(r)=\ell(\rho)<\ell(\rho')=\ell(r')$.
\pn
Suppose that $K=C(\eps_1, \ldots, \eps_\nu),\ \eps_i=\pm 1$. Let us show that
$\nu \geq n$. We have
$K=S(x)$ where
$x=[\eps_1, \ldots, \eps_\nu]$
and $x=\Frac{\alpha}{\beta+k\alpha}$ or $x=\Frac{\alpha}{\gamma+l\alpha}$ where $k,l \in \ZZ$.
We have $\ell(x)=\nu$ if $(\eps_1, \ldots, \eps_\nu)$ is 1-regular and $\nu\geq \ell(x)+3$ otherwise.
\bn
\item[--] If $k=2p>0$ then
we have $x=(MP)^p r$ so $\ell(x)=\ell(r)+3p>\ell(r)$.
\item[--] If $k=2p+1>0$ then
$x=(MP)^p M (\Frac 1r )$ so $\ell(x)=\ell(1/r)+3p+2=\ell(r')+3p+3>\ell(r')$.
\item[--] If $k=-(2p+1)<0$ then $-x=(MP)^p (r')$ so $\ell(x)=\ell(-x)=\ell(r')+3p>\ell(r')$.
\item[--] If $k=-2p>0$ then
$-x=(MP)^{p-1} M (\Frac 1{r'})$ so $\ell(x)=\ell(1/r')+3p-1=\ell(r)+3p>\ell(r)$.
\en
If $x=\Frac{\alpha}{\gamma+l\alpha}$, we obtain the same relations.
We deduce that $\nu\geq\min(\ell(r),\ell(r'))=n$ and the second point.
We deduce that $\nu=n$ if and only if
$x=\Frac\alpha\beta$ or $x=\Frac\alpha\gamma$. In this case, we get
the third point using Proposition \ref{gg} and the uniqueness
of the 1-regular continued fraction expansion (Theorem \ref{th1}).
\EPf
\begin{algorithm}[Computing the minimal Chebyshev diagram]\label{minb}\\
Let $K=S(\Frac\alpha\beta)$, $\Frac\alpha\beta>1$.
First compute the 1-regular sequence
$\Frac\alpha\beta = [e_1, \ldots, e_n],\ e_i=\pm 1$.
$C(e_1, \ldots, e_n)$ is
the Conway normal form corresponding to the Chebyshev diagram $\cC(3,n+1)$.
\pn
If $n< \frac 32 N -1$, then $b=n+1$ is the smallest integer such
that has a Chebyshev diagram $x=T_3(t),\, y=T_b(t)$, from Proposition \ref{ell}.
\pn
If $n> \frac 32 N -1$,
let us consider the 1-regular continued fraction
$\Frac\alpha{\alpha-\beta} = [\eps_1, \ldots, \eps_{n'}]$. We have
$n'<\frac 32 N -1$, by Proposition \ref{gg}.
$C(-\eps_1, \ldots, -\eps_{n'})$ is
the Conway normal form corresponding to the Chebyshev diagram $\cC(3,n'+1)$
of  $K = S(-\Frac\alpha{\alpha-\beta})$.
This last diagram is minimal by Proposition \ref{ell}.
\end{algorithm}
\begin{remark}[Minimality condition]\label{minib}
First compute $G \in \GG$ such that $\Frac\alpha\beta= PGP(\infty)$.
Let $p=\deg_P(PGP)$ and $m=\deg_M(PGP)$. By remark \ref{rl}, we have
$n=p+2m$ and $N=p+m$. Consequently, the minimality condition $n<\frac 32 N -1$ is equivalent
to $p\geq m+3$.
\end{remark}
\begin{example}[Torus knots]\label{tok}
The Schubert fraction of the torus knot $\TT(2,2k+1)$ is $2k+1$.
We have $ PM(x) = x+2,$ and then $ (PM)^k (x) = x+2k,  (PM)^kP (x)= 2k+1+ \Frac 1x.$
So we get the continued fraction of length $3k+1$:
$ 2k+1 = (PM)^kP(\infty)$.
This shows that the torus knot $\TT(2, 2k+1)$
has a Chebyshev diagram $ \cC( 3, 3k+2).$
This is not a minimal diagram.
\pn
On the other hand, we get $(PM)^{k-1}P^2(\infty) = 2k$ so
$\Frac{2k+1}{2k} = P(PM)^{k-1}P^2(\infty) >1$.
This shows that the torus knot $\TT(2, 2k+1)$
has a Chebyshev diagram $ \cC( 3,3k+1).$
This diagram is minimal by Remark \ref{minib}.
We proved in \cite{KP3} that $\TT(2, 2k+1)$ is in fact a harmonic knot
parametrized by $x=T_3(t),\, y= T_{3k+2}(t), \, z= T_{3k+1}(t)$, that is
$\deg y + \deg z = 3 (2k+1)$.
\end{example}
\begin{example}[Twist knots]\label{twk}
The twist knot $\cT_{n}$ is defined by $\cT_n = S(n+\Frac{1}{2})$.
\pn
From $P^3( x) = \Frac{3x+2}{2x+1}$, we get the continued fraction of length $3k+3$:
$\Frac {4k+3}2= (PM)^kP^3 (\infty).$
This shows that the twist knot $\cT_{2k+1}$
has a  Chebyshev diagram $\cC ( 3, 3k+4),$ which is minimal by Remark \ref{minib} .
\pn
We also deduce that $P(PM)^{k-1}P^3 (\infty) = P \left(\frac 12 (4k-1)\right ) = \Frac{4k+1}{4k-1}$.
This shows that the twist knot $\cT_{2k}$ has a  minimal Chebyshev diagram $\cC (3,{3k+2}) .$
\pn
We shall see that these knots are not harmonic knots for $a=3$
and we will give explicit bounds for their polynomial parametrizations.
\end{example}
\begin{example}[Generalized stevedore knots]\label{sk3}
The generalized stevedore knot ${\cal S}_k$ is defined by ${\cal S}_k = S(2k+2+\Frac{1}{2k})$.
We have
$$
(MP)^k = \cmatrix{1&0\cr 2k&1}, \, (PM)^k =\cmatrix{1&2k\cr 0&1}
$$
so $2k+2+\Frac{1}{2k} = (PM)^{k+1}(MP)^k (\infty)$. This shows
that the stevedore knot ${\cal S}_k$ has a Chebyshev diagram $\cC(3,6k+4)$.
It is not minimal and
we see, using Remark \ref{minib}, that the
knot ${\cal S}_k$  also has a minimal Chebyshev diagram $\cC(3,6k+2)$.
Moreover, using Proposition \ref{gg}, we get
$$
\Frac{(k+1)^2}{(k+1)^2-2k} = P^2 (MP)^k (PM)^{k-1} P^2 (\infty).
$$
\end{example}
\section{Polynomial parametrization of rational knots}\label{param}
\begin{definition}
Let ${\cal D}(K)$ be a diagram of a knot having crossing points corresponding to
the parameters  $ t_1, \ldots, t_{2m} $.
The Gauss sequence of ${\cal D}(K)$ is defined by
$ g_k = 1 $ if $t_k$ corresponds to an overpass, and
$ g_k = -1$ if $t_k$ corresponds to an underpass.
\end{definition}
\begin{theorem}\label{gauss3}
Let $K$ be a two-bridge knot of crossing number $N.$
Let $x= T_3(t),\ y =T_b(t)$ be the minimal Chebyshev diagram of $K$.
Let $c$ denote the number of sign changes in the corresponding Gauss sequence.
Then we have
$$ b+c=3N.$$
\end{theorem}
\Pf Let $s$ be the number of sign changes in the Conway normal form of $K.$
By Proposition \ref{bireg} we have $ N= b-1-s.$
From this we deduce that our condition is equivalent to $3s+c = 2b-3.$
Let us prove this assertion by induction on $s.$
If $s=0$ then the diagram of $K$ is alternating,
and we deduce $c= 2(b-1) -1=2b-3.$

Let $C(e_1,e_2, \ldots , e_{b-1}) $ be the Conway normal form of $K.$ We may suppose
$e_1=1$. We shall denote by $M_1, \ldots, M_{b-1}$ the crossing points of the diagram, and
by $ x_1< x_2 < \cdots < x_{b-1} $ their abscissae.
Let $e_k$ be the first negative coefficient in this form.
By the 1-regularity of the sequence we get $e_{k+1}<0,$ and $ 3 \le k \le b-1.$

Let us consider the knot $K'$ defined by its Conway normal form
$$ K' = C(e_1,e_2, \ldots , e_{k-1} , -e_k , -e_{k+1}, \ldots , - e_{b-1}).$$
We see that  the number of sign changes in the Conway sequence of $K'$ is
$s'= s-1.$ By induction, we get for the knot $K'$:  $3s'+c'= 2b-3.$

The  plane curve $x=T_3(t), \  y= T_b(t) $ is the union of three
arcs where $x(t)$ is monotonic.
Let $ \Gamma$ be one of these arcs.
$\Gamma$ contains (at least) one point $M_k$ or $M_{k+1}.$
\psfrag{g1}{\large$\mathbf{\Gamma_1}$}%
\psfrag{g2}{\large$\mathbf{\Gamma_2}$}%
\psfrag{g3}{\large$\mathbf{\Gamma_3}$}%
\psfrag{k1}{$M_{k-1}$}%
\psfrag{k2}{$M_k$}%
\psfrag{k3}{$M_{k+1}$}%
\begin{figure}[th]
\begin{center}
\begin{tabular}{ccc}
{\scalebox{.7}{\includegraphics{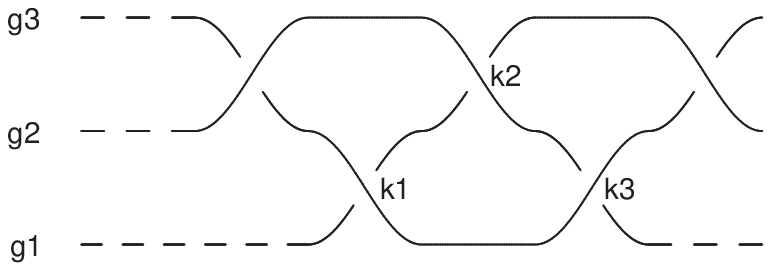}}}&&
{\scalebox{.7}{\includegraphics{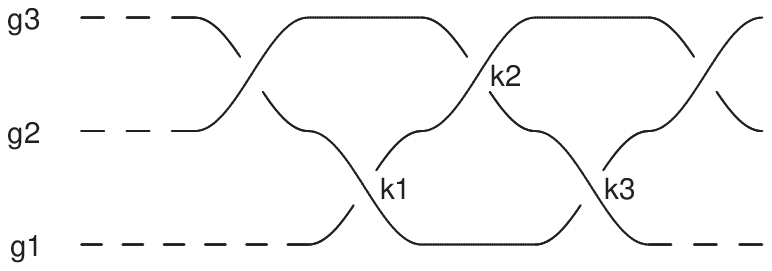}}}\\
$K$&\hspace{2cm}&$K'$
\end{tabular}
\end{center}
\caption{The modification of Gauss sequences}
\label{gaussf}
\end{figure}
Let $j$ be the first integer in $\{ k, k+1 \}$ such that $M_j$ is on $\Gamma$,
and let $j_- < j$ be the greatest integer such that $ M_{j_-} \in \Gamma .$
In figure \ref{gaussf}, we have
for $\Gamma_1$: $j=k, j^- = k-1$,
for $\Gamma_2$: $j=k, j^- = k-2$,
for $\Gamma_3$: $j=k+1, j^- = k-1.$

On each arc $\Gamma$, there is a sign change in the Gauss sequence iff the corresponding Conway
signs are equal.
Then, since the Conway signs $ s(M_{j_-} )$ and $s(M_j)$ are different, we see that
the corresponding Gauss signs are equal.
Now, consider the modifications in the Gauss sequences when we transform
$K$ into $K'$.
Since the the Conway signs $s( M_h) , \, h\ge k $ are changed,
we see that we get one more sign change on  every arc $\Gamma$.
Thus the number of sign changes in the Gauss sequence of $K'$ is
$c' = c+3.$
We get $3s+c= 3(s'+1) + c'-3= 3s'+c' = 2b-3,$ which completes our induction proof.
\EPf
\begin{corollary}\label{degc3}
Let $K$ be a 2-bridge knot with crossing number $N$. Then there exist $b,c$, $b+c=3N$,
and an  polynomial $C$ of degree $c$ such that the knot
$x=T_3(t), \, y=T_b(t), \, z=C(t)$ is isotopic to $K$.

If $K$ is amphicheiral, then  $b$ is odd, and the polynomial $C(t)$
can be chosen odd.
\end{corollary}
\Pf Let $b=n+1$ be the smallest integer  such that $K$ has
a Chebyshev diagram $x=T_3(t), \, y= T_b(t)$.
By our theorem \ref{gauss3}, the Gauss sequence $(g(t_1), \ldots, g(t_{2n}))$
of this diagram has $c=3N-b$ sign changes. We choose $C$ such that
$C(t_i)g(t_i)>0$ and we can realize it by choosing the roots of
$C$ in $]t_i,t_{i+1}[$ when $g(t_i)g(t_{i+1})<0$.

If $K$ is amphicheiral, then  $b$ is odd and the Conway form is palindromic
by Proposition \ref{palin}.
Then our Chebyshev diagram is symmetrical about the origin. We see that the Gauss
sequence is odd: $g(t_h) = -g(-t_h)$.
This implies that the polynomial $C(t)$ is odd when we choose its roots to be, for example,
$\frac 12(t_{i}+t_{i+1})$ where $g(t_i)g(t_{i+1})<0$.
\EPf
\begin{remark}
When $K$ if amphicheiral, it can be parametrized by three odd polynomials. In this case
the central symmetry $(x,y,z) \mapsto (-x,-y,-z)$ reverses the orientations of both $K$ and $\SS^3$.
This gives a simple proof of a famous theorem
of Hartley and Kawauchi:
{\em every amphicheiral rational knot is strongly negative amphicheiral}
(\cite{HK,Kaw}).
\end{remark}
Theorem \ref{gauss3} provides an effective polynomial parametrization $(x(t),y(t),z(t)$ with
$\deg x = 3$, $\deg y + \deg y = 3N$. We conjecture:
\begin{conjecture}
Let $K$ be a rational knot of crossing number $N$. Let $(x(t),y(t),z(t))$ be a polynomial
parametrization of $K$ with $\deg x =3$, then we have $\deg y + \deg z \geq 3N$.
\end{conjecture}
We shall give several examples of polynomial parametrizations of rational knots
with Chebyshev diagrams
$\cC(3,b)$.
\subsubsection*{Parametrizations of the torus knots}
The torus knot $\TT(2,2n+1) = S(2n+1)$ has a minimal Chebyshev diagram
$\cC(3,3k+1)$. It can be parametrized by $z=C(t)$ where
$\deg (C) = 3k+2$. Actually we proved in \cite{KP3} that $C = T_{3k+2}$ is
convenient.
\subsubsection*{Parametrizations of the twist knots}
The twist knot $\cT_m= S(m + \Frac 12 ) $ has crossing number $m+2$.
We shall see that $\cT _m$  is not
a harmonic knot $ \H(3,b,c)$  because $2^2 \not \equiv \pm 1 \Mod{2m+1}$ except when
$m=2$ (the figure-eight knot) or $m=1$ (trefoil).
From example \ref{twk}, we know that:
\bi
\item $\cT_{2k+1}$ can be parametrized by
$x=T_3(t), \, y=T_{3k+4}, \, z=C(t)$ where $\deg(C)=3k+5$.
\item $\cT_{2k}$ can be parametrized by
$x=T_3(t), \, y=T_{3k+2}, \, z=C(t)$ where $\deg(C)=3k+4.$
\ei
\begin{example}[The 3-twist knot]
$\cT_3$ is the 3-twist knot $\overline{5}_2$. It is the harmonic knot $\H(4,5,7)$
(see \cite{KP3}).
It can also be parametrized by
$$x=T_3(t), \ y=T_7(t),\ z= t \left( 4\,t+3 \right)  \left( 3\,t+1 \right)  \left( 6\,t-5 \right)
 \left( 12\,{t}^{2}-11 \right)  \left( 2\,{t}^{2}-1 \right)$$
\begin{figure}[th]
\begin{center}
\begin{tabular}{cc}
{\scalebox{.8}{\includegraphics{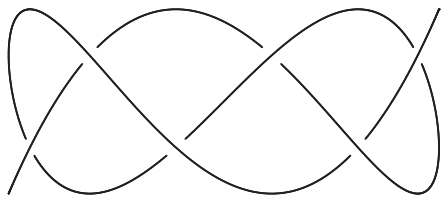}}}
\end{tabular}
\end{center}
\caption{Diagram of the 3-twist knot $\overline{5}_2=S(\Frac 72)$}
\label{K52}
\end{figure}
\end{example}
\subsubsection*{Parametrizations of the generalized stevedore knots}
The stevedore knot ${\cal S}_m= S (2m+2+ \Frac 1{2m} )$  can be represented by
$ x=T_3(t), \ y= T_{6m+2}(t), \ z = C(t) $ where $C(t)$ is a polynomial
of degree $6m+4$.
\pn
\Pf
This is a consequence of  \ref{sk3} and Corollary \ref{degc3}.
\EPf
\begin{example}[The knot $6_1={\cal S}_1$]\label{6_1-3}
In the example \ref{9297}, we get $\ell(\Frac 92)=9,\, \ell(\Frac 97)=7$. $b=8$ is the minimal value for which
$x=T_3(t), \, y=T_8(t)$ is a Chebyshev diagram for $\overline{6}_1$.
The Gauss sequence associated to the Conway form
$\overline{6}_1 = C(-1,-1,-1,1,1,1,1)$ has exactly 10 sign changes. It is precisely
$$[1, -1, -1, 1, -1, 1, -1, -1, 1, -1, 1, 1, -1, 1].$$
We can check that
$$
x=T_3(t),\, y=T_8(t),\, z=
 \left( 8\,t+7 \right)  \left( 5\,t-4 \right)  \left( 15\,{t}^{2}-14
 \right)  \left( 2\,{t}^{2}-1 \right)  \left( 3\,{t}^{2}-1 \right)
 \left( 15\,{t}^{2}-1 \right)
$$
is a parametrization of $\overline{6}_1$ of degree $(3,8,10)$.  In \cite{KP3} we gave
the Chebyshev parametrization $6_1= \cC(3,8,10,\frac{1}{100})$.
We shall see that $\overline{6}_1= S(\Frac 92)$ is not a harmonic knot $ \H(3,b,c)$
because $2^2 \not \equiv \pm 1 \Mod 9$.
\begin{figure}[th]
\begin{center}
{\scalebox{.7}{\includegraphics{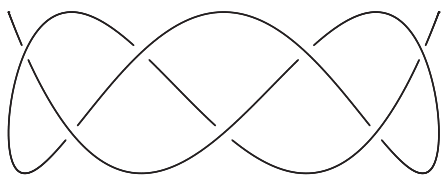}}}
\end{center}
\caption{The knot $6_1$}
\label{61_3}
\end{figure}
\end{example}
\section{The harmonic knots $\mathbf{\H(3,b,c)}$}\label{harmonic}
In this paragraph we shall study Chebyshev knots with $ \phi= 0.$
Comstock (1897) found the number of crossing points of the
harmonic  curve parametrized by $x=T_a(t), y=T_b(t), z=T_c(t)$.
In particular, he proved that this curve is non-singular if and
only if $ a,b,c$ are pairwise coprime integers (\cite{Com}).
Such harmonic curves will be named harmonic knots  $\H(a, b,c)$ following the
original denomination (\cite{Com}).
These are not the harmonic knots defined by Trautwein (\cite{Tr}),
which are now referred to as Fourier knots (cf \cite{Cr}).
\pn
We shall need the following result proved in \cite{KP3}
\begin{proposition}
Let $a $ and $b$ be coprime integers.  The $\frac 12 (a-1)(b-1)$
double points of the Chebyshev curve
$x= T_a(t), y= T_b(t)$ are obtained for the parameter pairs
$$
t= \cos \bigl( \Frac ka + \Frac hb \bigr) \pi,  \
s = \cos \bigl( \Frac ka - \Frac hb \bigr) \pi ,
$$
where $h,k$ are positive integers such that
$ \Frac ka + \Frac hb < 1 .$
\end{proposition}
Using the symmetries of Chebyshev polynomials, we see
that this set of parameters is symmetrical about the origin.
We shall need the following result proved in \cite{KP3}.
We will write  $ x \sim y$ when  $\sign x = \sign y .$
\begin{lemma}\label{sign}
Let $\H(a,b,c)$ be the harmonic knot: $x=T_a(t),\, y=T_b(t),\, z=T_c(t)$.
A crossing point of parameter
$
t = \cos \left ({k \over a } + {h \over b } \right ) \pi,  \
$
is a right twist if and only if
$$D = \Bigl( z(t)-z(s) \Bigr) x'(t) y'(t) >0$$
where
$$
z(t)-z(s) = T_c(t)-T_c(s) =
-2 \sin \Bigl( \Frac {ch}{b} \pi \Bigr) \sin \Bigl( \Frac {ck}a \pi \Bigr).
$$
and
$$
x'(t) y'(t) \sim (-1)^{h+k} \sin \Bigl( \Frac {ah}b \pi \Bigr)
\sin \Bigl( \Frac {bk}a \pi \Bigr) .
$$
\end{lemma}
From this lemma we immediately deduce
\begin{corollary}\label{cprime}
Let $a,b,c$ be coprime integers. Suppose that the integer $c'$ verifies
$ c' \equiv c  \Mod{2a} $ and $ c' \equiv -c \Mod{2b} .$
Then the knot $\H(a,b,c')$ is the mirror image of $\H( a,b,c).$
\end{corollary}
\Pf
At each crossing point  we have
$
T_{c'}(t) - T_{c'}(s) = - \Bigl( T_c(t) - T_c (s)  \Bigr).
$
\EPf
\begin{corollary}
Let $a,b,c$ be coprime integers. Suppose that the integer $c$ is of
the form $c= \lambda a + \mu b$ with $\lambda, \mu >0$.
Then there exists $c'< c $
such that
$\H( a,b,c)= \overline{\H} (a,b,c') $
\end{corollary}
\Pf
Let $c'=\abs{\lambda a - \mu b}.$
The result follows immediately from corollary \ref{cprime}
\EPf
\pn
In a recent paper, G. and J. Freudenburg  have proved the following stronger result.
{\em There is a polynomial automorphism $ \Phi$
of $ \RR^3$ such that $ \Phi (\H(a,b,c)) = \H(a,b,c').$}
They also conjectured that the knots
$\H(a,b,c), \  a<b<c,  \  c \neq \lambda a + \mu b , \  \lambda, \mu >0 $
are  different knots (\cite{FF}, Conjecture 6.2).
\pn
When $a\leq 2$, it is easy to see that the harmonic knots $\H(a,b,c)$ are trivial knots.
\pn
The following result is the main step in the classification of the
harmonic knots $\H( 3,b,c)$.
\begin{theorem}\label{h3}
Let $b=3n+1, \ c= 2b-3 \lambda , \, (\lambda, b)=1 .$
The Schubert fraction of the knot $\H (3, b,c )$
is
$$ \Frac {\alpha}{\beta} =
[e_{1}, e_{2}, \ldots, e_{3n}],
\hbox{ where } e_k = \sign {\sin k \theta} \hbox{ and } \theta = \Frac {\lambda}b  \pi. $$
If $0< \lambda < \Frac {b}{2}$, its crossing number is
$N= b -  \lambda = \Frac {b+c}3,$
and we  have $\beta^2 \equiv \pm 1 \Mod \alpha.$
\end{theorem}
\Pf Will be given in section {\bf \ref{proofs}}, p. \pageref{proofs}.
\begin{corollary}\label{hh3}
The knots $\H(3,b,c)$ where $\Frac c2 <b<2c, \  b  \equiv 1 \Mod 3, \  c \equiv 2 \Mod 3$
are different knots (even up to mirroring). Their crossing number is given by $b+c=3N.$
\end{corollary}
\Pf
Let $K=\H(3,b,c)$ and $\Frac \alpha \beta>1$
be its 1-regular Schubert fraction given by Theorem \ref{h3}.
From Prop \ref{ell}, $\min(b,c)$ is
the minimum length of a Chebyshev diagram of  $K$ and
$\max(b,c)=3N-\min(b,c)$. The pair $(b,c)$ is uniquely determined.
\EPf
\pn
The following result gives the classification of harmonic knots $\H(3,b,c)$.
\begin{theorem}\label{h3bc}\\
Let $K= \H(3,b,c).$
There exists  a unique  pair $(b',c')$ such that (up to mirror symmetry)
$$
K = \H(3,b',c'),   \  b'<c'<2b', \  b'+c' \equiv 0 \Mod 3.
$$
The crossing number of $K$ is $\frac 13 (b'+c'),$
its fractions $ \Frac \alpha \beta $ are such that
$ \beta^2 \equiv \pm 1 \Mod \alpha.$
Furthermore, there is an algorithm to find the pair $(b',c').$
\end{theorem}
\Pf
Let $K= \H(a,b,c)$
We will show that if the pair $(b,c)$ does not satisfy the condition of the theorem,
then it is possible to reduce it.

If $c<b$ we consider $\H(3,c,b)=\overline \H(3,b,c)$.

If $b \equiv c \Mod 3,$ we have $ c= b + 3 \mu, \  \mu>0.$
Let $c'= \abs{b-3\mu}$.
We have $c'\equiv \pm c \Mod {2b}$ and  $c' \equiv\mp c \Mod 6.$
By Lemma \ref{cprime}, we see that $ K= \overline \H(3,b,c')$ and we get a smaller
pair.

If $b \not\equiv c \Mod 3$ and $c>2b,$ we have $c=2b + 3\mu, \ \mu>0.$
Let $c'= \abs{2b-3 \mu}$.
Similarly, we get  $ K= \overline\H(3,b,c')$.
This completes the proof of  existence.
This uniqueness is a direct consequence of  Corollary \ref{hh3}.
\EPf
\begin{remark}
In \cite{FF}, it is proved (see Proposition 4.2) that $\H(3,b,c)$ and $\H(3,b',c')$ are
algebraically equivalent.
\end{remark}
\begin{remark}
Theorem \ref{h3bc} gives a positive answer to the Freudenburg conjecture for $a=3.$
\end{remark}
\begin{remark}
Harmonic knots have Chebyshev parametrizations with the same degrees as the
parametrizations we gave in section {\bf\ref{param}} (Theorem \ref{gauss3}).
\end{remark}
\begin{example}
We get $\H(3,31,43) = \H(3,19,31) = \H(3,7,19) = \H(3,5,7)$. The crossing number of
this knot is $4=\frac 13 (5+7)$. With $b=5$ and $c=2b-3$ we get the Conway normal form
$C(\sign{\sin \Frac{\pi}5}, \sign{\sin \Frac{2\pi}5},\sign{\sin \Frac{3\pi}5},\sign{4\sin \Frac{\pi}5})
=C(1,1,1,1)$. Its Schubert fraction is $\Frac 53 = [1,1,1,1]$. It is the knot $4_1$.
\end{example}
\subsection*{Examples}
As  applications of Proposition \ref{bireg},
 let us deduce the following results (already   in \cite{KP3}).
\begin{corollary}
The harmonic knot $ \H (3, 3n+2, 3n+1) $ is the torus knot $ \TT(2, 2n+1).$
\end{corollary}
\Pf
The harmonic knot $ K= \H (3, 3n+1, 3n+2 )$ is obtained for
$ b= 3n+1,$ $c= 2b -3\lambda, \  \lambda= n ,$
 $ \theta = \Frac n {3n+1} \pi .$
If $j=1,2,$ or $ 3,$ and $ k=0,\ldots , n-1$
we have
$(3k+j) \theta = k \pi + \Frac {jk-n}{3n+1} ,$
hence $ \sign{ \sin (3k+j) \theta } = (-1)^k $, so that the Schubert fraction of $K$ is
$$
[1,1,1, -1,-1,-1, \ldots , (-1)^{n+1}, (-1)^{n+1} , (-1)^{n+1} ] = \Frac{2n+1} {2n}
\thickapprox -(2n+1) .
$$
We see that $K$ is the mirror image of $ \TT(2,2n+1),$ which completes the proof.
\EPf
\begin{figure}[th]
\begin{center}
\begin{tabular}{ccc}
{\scalebox{.7}{\includegraphics{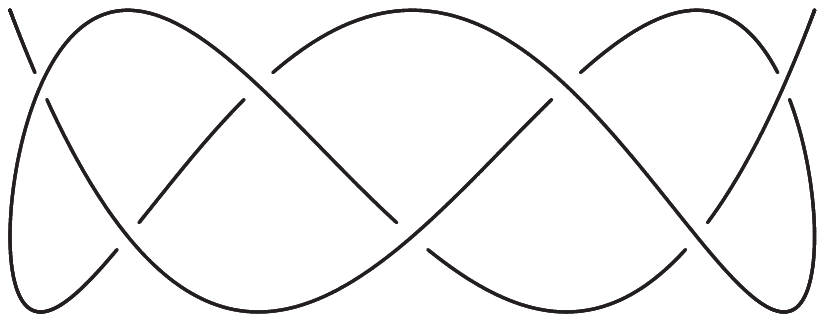}}}&&
{\scalebox{.7}{\includegraphics{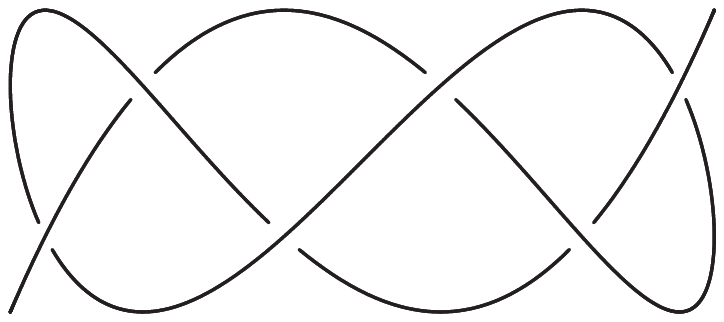}}}\\
$\H(3,8,7)$&\hspace{2cm}&$\H(3,7,8)$
\end{tabular}
\end{center}
\caption{The torus knot $\TT(2,5)=5_1$ and its mirror image}
\label{t25}
\end{figure}
\pn
It is possible to parameterize the knot $\TT(2, 2n+1)$ by polynomials of the same degrees
and an alternating diagram (\cite{KP2}).
However, our Chebyshev parametrizations are easier to visualize.
We conjecture that these degrees are minimal
(see also \cite{RS,KP1}).
\begin{corollary}
The harmonic knot $\H(3,b,2b-3)$ ($b\not \equiv 0 \Mod 3$)  has
crossing number $b-1$. The Chebyshev diagram of the projection on the
$xy$-plane is alternating.
\end{corollary}
\Pf
For this knot we have $\lambda=1, \  \theta = \Frac \pi b .$
The Conway normal form of the projection on the
$xy$-plane is $(1,1,\ldots, 1)$.
The Schubert fraction is the continued fraction of length $b-1$:
$[1,1, \ldots ,1] = \Frac{F_{b}}{F_{b-1}}$ where $F_n$ are the Fibonacci numbers
($F_0=0, F_1=1, \ldots$).
\EPf
\begin{remark}
J. C. Turner named these  knots  Fibonacci knots (\cite{Tu}).
In \cite{KP3}, we showed the more general result: the projection
of $\H(a,b,ab-a-b)$ on the $xy$-plane is alternating.
In \cite{KP5}, we have
studied Fibonacci knots and generalized Fibonacci knots and showed
that most of them are not Lissajous knots.
\end{remark}
The two previous examples describe infinite families of harmonic knots. They have
a Schubert fraction $\Frac\alpha\beta$ with $\beta^2 \equiv  1 \Mod\alpha$ (torus knots) or
with $\beta^2 \equiv - 1\Mod\alpha$ (Fibonacci knots with odd $b$). There
is also an infinite number of
two-bridge knots with $\beta^2= \pm 1 \Mod\alpha$ that are not harmonic.
\begin{proposition}
The knots (or links)
$K_n  = \cC(1,1,\underbrace{-1,\ldots,-1}_{n+2},1,1), \ n>1$,
are not harmonic knots $\H(3,b,c)$. Their crossing number is $n+4$ and their
Schubert fraction $\Frac{\alpha_n}{\beta_n}$ satisfy $\beta_n^2 \equiv  (-1)^{n+1} \Mod{\alpha_n}$.
\end{proposition}
\Pf
We have
$
\Frac{\alpha_n}{\beta_n}=P M P^n M P(\infty).$
Using the fact that $P^n = \cmatrix{F_{n+1}&F_n\cr F_n&F_{n-1}}$ we deduce
that
$$
P M P^n M P = \cmatrix{5F_{n+1}&F_{n+1}+F_{n-1}\cr F_{n+1}+F_{n-1}&F_{n-1}},
$$
that is $\alpha_n =5F_{n+1}, \, \beta_n =F_{n+1}+F_{n-1}$.
Taking  determinants, we get $\beta_n^2 \equiv (-1)^{n+1} \Mod{\alpha_n}$.
Since $ n+2 \ge 4,$ it cannot be of the form
  $[\sign{\sin \theta},\sign{\sin 2\theta}, \ldots, \sign{\sin k\theta}]$.

If $n\equiv 2\Mod 3$, $K_n$ is a two-component link.

If $n\equiv 1\Mod 6$ or $n\equiv 3\Mod 6$, we have
$\beta_n^2\equiv 1\Mod{\alpha_n}$.

If $n\equiv 0\Mod 6$ or $n\equiv 4\Mod 6$, $K_n$ is amphicheiral.
\EPf

\section{Proof of theorem \ref{h3} }\label{proofs}
We study here the diagram of $\H(3,b,c)$ where $b=3n+1$ and
$c=2b-3\lambda$. The crossing points of the plane projection of
$\H(3,b,c)$ are obtained for pairs of values $(t,s)$
where
$t= \cos \bigl( \Frac m{3b} \pi \bigr) , \ s = \cos \bigl( \Frac {m'}{3b} \pi \bigr).$
For $k = 0, \ldots, n-1$, let us consider
\bi
\item $A_{k}$  obtained for $m=3 k +1, \  m'= 2b-m.$
\item $B_{k}$  obtained for $ m =   3 k +2, \ m'= 2b + m$.
\item $C_{k}$  obtained for $ m = 2b - 3k - 3, \  m'= 4b-m$.
\ei
Then we have
\bi
\item $x(A_{k}) = \cos \bigl( \Frac {3 k +1}b \pi \bigr)$, \
$y(A_k) = \frac 12 (-1)^k$.
\item $x(B_{k} )= \cos \bigl( \Frac {3 k +2}b \pi \bigr)$, \
$y(B_k) = \frac 12 (-1)^{k+1}$.
\item
$x(C_{k}) = \cos \bigl( \Frac {3 k +3}b \pi \bigr)$, \
$y(C_k) = \frac 12 (-1)^{k}$.
\ei
\psfrag{a0}{\small $A_0$}\psfrag{b0}{\small $B_0$}\psfrag{c0}{\small $C_0$}%
\psfrag{a1}{\small $A_1$}%
\psfrag{an}{\small $A_{n-1}$}\psfrag{bn}{\small $B_{n-1}$}\psfrag{cn}{\small $C_{n-1}$}%
\begin{figure}[th]
\begin{center}
{\scalebox{.7}{\includegraphics{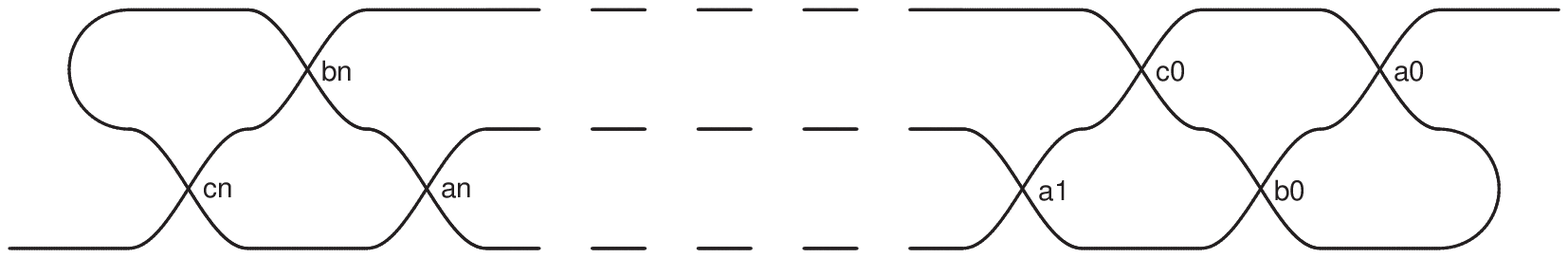}}}
\end{center}
\caption{$\H(3,3n+1,c)$, $n$ even}
\label{dh3}
\end{figure}
Hence our $3n$ points satisfy
$$x(A_{k-1}) > x(B_{k-1})
>x(C_{k-1}) >  x( A_{k}) > x(B_{k}) > x ( C_{k} ),\ k=1, \ldots , n-1 .$$
Using the identity $T'_a( \cos \tau ) =  a  \Frac{ \sin a \tau}{\sin \tau },$
we get
$
x'(t)y'(t)  \sim \sin \bigl( \Frac mb \pi \bigr)  \sin \bigl( \Frac m3 \pi \bigr) .
$
We obtain
\bi
\item[]
for $A_k$:
$
\begin{array}[t]{rcl}
{x'(t)y'(t)}&\sim& {\sin ( \Frac {3 k +1}b \pi ) \sin ( \Frac {3 k
+1}3 \pi )}  \sim (-1)^{k}.
\end{array}$
\item[] for $B_k$:
$
\begin{array}[t]{rcl}
x'(t)y'(t)&\sim& {\sin \bigl( \Frac {3 k +2}b \pi \bigr) \sin
\bigl( \Frac{3k+2}3 \pi \bigr)}
 \sim (-1)^{k}.
\end{array}$
\item[] for $C_k$:
$
\begin{array}[t]{rcl}
{x'(t)y'(t)}&\sim&
{\sin ( \Frac {2b -3 k -3}b \pi ) \sin (\Frac {2b-3
k -3}3 \pi )}\\
&\sim&
- {\sin (\Frac{3k+3}b \pi ) \sin ( -\Frac {3k+1}3  \pi ) } \sim
(-1)^{k}.
\end{array}$
\ei
The following identity will be useful in computing the sign of $z(t)-z(s).$
$$
T_c(t) - T_c(s) = 2 \sin \Bigl( \Frac{c}{6b}(m'-m) \pi  \Bigr)
\sin \Bigl( \Frac{c}{6b}(m+m') \pi   \Bigr).
$$
We have, with $c=2b -3\lambda$, $\theta= \Frac{\lambda}b \pi ,$ (and $b=3n+1$ ),
\bi
\item[]
for $A_k$:
$ z(t)-z(s) =
-2 \sin c \Frac{\pi}3  \sin \Bigl( c \Frac{m-b}{3b} \pi \Bigr).$
But
\[
\sin  c \Frac {\pi}3 = \sin \Bigl( \Frac {6n+2-3 \lambda}3 \pi \Bigr)
= (-1)^{\lambda}  \sin \Frac {2 \pi }3,\label{c}
\]
and
$$ \sin \Bigl( c \Frac{b-m}{3b} \pi \Bigr)=
\sin \Bigl( (2- \Frac {3 \lambda}b )  \,  \Frac {b-m}3 \pi \Bigr) =
\sin \bigl( \Frac {\lambda}b (m-b) \pi \bigr) =
(-1)^{\lambda} \sin ( 3k+1) \theta .$$
We deduce that $z(t) - z(s) \sim \sin (3k+1) \theta$.
Finally, we obtain
$$ \sign{D(A_k)}= (-1)^k\sign{\sin (3k+1) \theta}.$$
\item
for $B_k$:
$
z(t)-z(s) =
2 \sin  c \Frac{\pi}{3}  \sin \bigl( \Frac cb \, . \,  \Frac{b+m}{3} \pi\bigr).
$
We have
\[
\sin \bigl( \Frac cb \cdot \Frac {b+m}3 \pi \bigr) &=&
\sin \bigl( ( 2-\Frac {3\lambda}b ) \Frac {b+m}3 \pi \bigr) \nonumber
\\
&=&
-\sin \bigl( \Frac {\lambda}b (b+m) \pi \bigr) =
(-1)^{\lambda+1} \sin  (3k+2) \theta. \nonumber
\]
Then, using Equation \ref{c}, we get $z(t) - z(s) \sim - \sin( 3k+2) \theta$,
and finally
$$
\sign{D(B_k)} = (-1)^{k+1} \sign{\sin(3k+2) \theta}.
$$
\item
for $C_k$:
$
\begin{array}[t]{rcl}
z(t)-z(s) &\sim&
\sin  \Frac {2 c }3 \pi  \sin \bigl( \Frac cb ( k+1)  \pi \bigr)\\
&\sim&  \sin   \Frac {4 \pi } 3  \sin \bigl(
(2- \Frac{3\lambda}b )( k+1) \pi \bigr) \sim     \sin (3k+3) \theta.
\end{array}$\\
We obtain
$$
\sign{ D(C_k)} = (-1)^k \sign{ \sin (3k+3) \theta } .
$$
\ei
These results give the Conway normal form.
If $n$ is odd, the Conway's signs of our points are
$$
\begin{array}{rcccl}
s(A_k) &\sim& (-1)^k D( A_k ) &\sim& { \sin (3k+1) \theta },\\
s(B_k)&\sim& (-1)^{k+1} D(B_k) &\sim& { \sin (3k+2) \theta }, \\
s(C_k) &\sim& (-1)^k D(C_k) &\sim& { \sin (3k+3) \theta } .
\end{array}
$$
In this case our result follows,  since the fractions
$ [a_1, a_2, \ldots , a_{3n} ] $ and  $(-1)^{3n+1} [a_{3n}, \ldots , a_1] $ define
the same knot.
If $n$ is even, the Conway's signs are
the opposite signs, and we also get the Schubert fraction of our knot.

Since $ 0<\theta < \Frac {\pi}2 $, we see that
there are not two consecutive sign changes in our sequence.
We also see that the  first two terms are of the same sign, and so are the last two terms.
The Conway normal form is 1-regular and the total number of sign changes in this sequence is
$\lambda -1 $: the crossing number of our knot is then
$b- \lambda.$ Finally, we get $ \beta ^2 \equiv \pm 1 $ by Proposition \ref{palin}.
\EPf
\section{Conclusion}
We have given here a complete classification of harmonic knots $\H(3,b,c)$ by
computing explicitly their Schubert fraction.
We have shown that when $b<c<2b$ then $\H(3,b,c)$ has crossing number
$N=\frac 13 (b+c)$.
\pn
On the other hand we have shown that any rational knot of
crossing number $N$ admits a polynomial parametrization of degrees $(a,b,c)$ where
$a=3$, $N = \frac 13 (b+c)$ and $N<b<c<2N$.
\pn
This is the first algorithm giving explicit polynomial parametrizations for the
infinite family of rational knots.
We also conjecture that these degrees are minimal ($a=3$, $b+c=3N$).


\begin{thebibliography}{ZZ99}
%
\bibitem [BZ]{BZ}
G. Burde, H. Zieschang,
{\it Knots}, Walter de Gruyter, 2003
%
\bibitem [BDHZ] {BDHZ}
A. Boocher, J. Daigle, J. Hoste, W. Zheng,
{\it Sampling Lissajous and Fourier knots},
arXiv:0707.4210, (2007).
%
\bibitem[BHJS]{BHJS}
M. G. V. Bogle, J. E. Hearst, V. F .R. Jones, L. Stoilov, {\it
Lissajous knots}, Journal of Knot Theory and its Ramifications,
3(2): 121-140, (1994).
%
\bibitem [Com]{Com}
E. H. Comstock,
{\it The Real Singularities of Harmonic Curves of three Frequencies},
Trans. of the Wisconsin Academy of Sciences, Vol XI : 452-464, (1897).
%
\bibitem [Con]{Con}
J. H. Conway,
{\it An enumeration of knots and links, and some of their algebraic properties},
Computational Problems in Abstract Algebra (Proc. Conf., Oxford, 1967), 329--358 Pergamon,
Oxford (1970)
%
\bibitem [Cr]{Cr}
P. R. Cromwell, {\it Knots and links},
Cambridge University Press, Cambridge, 2004. xviii+328 pp.
%
\bibitem[Fi]{Fi}
G. Fischer, {\it Plane Algebraic Curves}, A.M.S. Student Mathematical
Library Vol 15, 2001.
%
\bibitem[FF]{FF}
G. Freudenburg, J. Freudenburg,
{\it Curves defined by Chebyshev polynomials}, 19 p., (2009),
{\tt arXiv:0902.3440}
%
\bibitem[HK]{HK}
R. Hartley, A. Kawauchi, {\it Polynomials of amphicheiral knots},
Math. Ann. {\bf 243 (1)}: 63-70 (1979)
%
%
\bibitem[HZ]{HZ}
J. Hoste, L. Zirbel, {\it Lissajous knots and knots with Lissajous
  projections}, (2006),
{\tt arXiv:math/0605632}.
To appear in Kobe Journal of mathematics, vol 24, n$^{\rm o}2$
%
\bibitem[JP]{JP}
V. F. R. Jones, J. Przytycki, {\it Lissajous  knots and billiard knots,}
Banach Center Publications, 42:145-163, (1998).
%
\bibitem[Kaw]{Kaw}
A. Kawauchi, editor, {\it A Survey of Knot Theory}, Birh{\"a}user, 1996.
%
\bibitem[KP1]{KP1}
P. -V. Koseleff, D. Pecker, {\it On polynomial Torus Knots},
Journal of Knot Theory and its Ramifications, Vol. {\bf 17 (12)}
(2008), 1525-1537.
%
\bibitem[KP2]{KP2}
P. -V. Koseleff, D. Pecker, {\it A construction of polynomial torus
 knots}, to appear in Journal of AAECC,
{\tt arXiv:0712.2408}.
%
\bibitem[KP3]{KP3}
P. -V. Koseleff, D. Pecker, {\it Chebyshev knots},  (2008), {\tt arXiv:0812.1089}.
%
\bibitem[KPR]{KPR}
P. -V. Koseleff, D. Pecker, F. Rouillier, {\it The first rational Chebyshev knots}, Conference MEGA 2009, Barcelona.
%
\bibitem[KP4]{KP4}
P. -V. Koseleff, D. Pecker, {\it Chebyshev diagrams for rational knots}, (2008),
{\tt arXiv:0906.4083}.
%
\bibitem[KP5]{KP5}
P. -V. Koseleff, D. Pecker, {\it On Fibonacci knots}, (2009),
{\tt arXiv:0908.0153}.
%
\bibitem[La1]{La1}
C. Lamm, {\it There are infinitely many Lissajous knots,}
Manuscripta Math., 93: 29-37, (1997).
%
\bibitem[La2]{La2}
C. Lamm, {\it Cylinder knots and symmetric unions
(Zylinder-knoten und symmetrische Vereinigungen),} Ph.D. Thesis,
Bonner Mathematische Schriften 321, Bonn, 1999.
%
%
%
%
\bibitem[Mu]{Mu}
K. Murasugi, {\it Knot Theory and its Applications}, Boston,
Birkh{\"a}user, 341p., 1996.
%
\bibitem[P1]{P1}
D. Pecker, {\it Simple constructions of algebraic curves with nodes},
Compositio Math.  87 (1993),  no. 1, 1--4.
%
\bibitem[RS]{RS}
A. Ranjan and R. Shukla, {\it On polynomial representation of
torus knots,} Journal of knot theory and its ramifications, Vol. 5
(2) (1996) 279-294.
%
%
\bibitem[Sh]{Sh}
A.R. Shastri, {\it Polynomial representation of knots},
 T{\^o}hoku Math. J. {\bf 44 } (1992), 11-17.
%
\bibitem[St]{St}
A. Stoimenow, {\it Generating functions, Fibonacci numbers and rational knots},
J. Algebra {\bf 310(2)} (2007), 491--525.
%
\bibitem[Tr]{Tr}
A. Trautwein, {\it Harmonic knots},
Ph.D. Thesis, University of Iowa, 1994
%
\bibitem[Tu]{Tu}
J.C. Turner, {\it On a class of knots with Fibonacci invariant numbers},
Fibonacci Quart. {\bf 24} (1986), n$^{\rm o}1$,  61-66.
%
\bibitem[Va]{Va}
V. A. Vassiliev, {\it Cohomology of knot spaces}, Theory of
singularities and its Applications, Advances Soviet Maths Vol 1
(1990)
\end{thebibliography}
\end{document}